\documentclass{article}
\usepackage{pdfpages}
\usepackage{graphicx}
\usepackage{amssymb}	
\usepackage{amsmath}
\begin{document}

\centerline{\large\bf Optimal stabilization of a cycle in nonlinear discrete systems}
\medskip

\centerline{\bf D.Dmitrishin, A.Khamitova, A.Korenovskyi and A. Stokolos} 
\bigskip

\textbf{1. Motivation.} 

 The problem of optimal influence on a chaotic regime is the most fundamental in non-linear dynamics (c.f.\cite{C,U} or \cite{M} for recent updates). The aim is either synchronizing of a chaotic motion or conversely, make a regular motion to be chaotic. Moreover, the admissible controls are only small ones which however totally change the regime of the motion. The solution to the problem of optimal stabilizing of an equilibrium point in a discrete autonomous system with the classical delayer feedback control (DFC) was given in \cite{DK,DH}. Earlier results in this direction had certain limitations (c.f. \cite{U}). It turns out that a machinery developed in \cite{DK,DH} still work for a cycle as well.

{\bf 2. Statement of the problem.} A central problem in the dynamic systems theory is to study a finite parameter family of mappings and to study a dependence of the qualitative properties of the generated systems on those parameters.

Namely, let us consider an open scalar non-linear discrete system
\begin{equation} \label{1} 
x_{n+1} =f\left(x_{n} \right),\, \, x_{n} \in \mathbb R^{1} ,\, \, n=1,\, \, 2,\, \, \ldots \, \, , %
\end{equation} 
with a non-stable cycle $\left(\eta _{1} , \eta _{2} \right)$ which does mean that  $\eta _{1} ,\eta _2 $ are distinct and that $\eta _{2} =f_{h} \left(\eta _{1} \right)$ while $\eta _{1} =f_{h} \left(\eta _2 \right).$ 

 It is assumed that a differentiable function  $f_h$ depends on a vector of parameters $h$ and that for every admissible set of these parameters the function is defined on a certain bounded interval and maps it into itself. Hence the cycle $\left(\eta _{1} ,\eta _2 \right)$ will dependent on these parameters as well as the cycle multiplier  
 $\mu =f_h^\prime\left(\eta _{1} \right)\cdot f_h^\prime(\eta _{2}).$ 
 
It is assumed that $\mu \in \left(-\mu ^{*} ,\, -1\right),\, \, \mu ^{*} >1$, and that for some $\mu \in \left(-\mu ^{*} ,\, -1\right)$ the phenomena of quasi-dynamical chaos is observed.
We would like to suppers the chaos by stabilization of the cycle for all admissible parameters  $h$  by the control of the form
\begin{equation} \label{2} 
u=-\sum _{j=1}^{N-1}\varepsilon _{j} \left(f_{h} \left(x_{n-2j+2} \right)-f_{h} \left(x_{n-2j} \right)\right), |\varepsilon _{j}|<1, j=1,\ldots, N-1, 
\end{equation} 
in a such way that the depth of the used prehistory  $N^{*} =2(N-1)$ is minimal.\\ 
 
Let’s indicate, that after the synchronization of the state ${x_n} = {x_{n - 2}}$   the control \eqref{2} resets, i.e. closed system takes the form, as in the absence of control. It means, that a cycles of the open and closed systems are the same.

{\bf 3. Auxiliary results.}
Let rewrite  \eqref{2} in the following form
\begin{equation} \label{3} 
x_{n+1} =\sum _{j=1}^{N}a_{j} f \left(x_{n-2j+2} \right) \, ,                                                             
\end{equation} 
where $a_{1} =1-\varepsilon _{1} , a_{j} =\varepsilon _{j-1} -\varepsilon _{j} , j=2,\, \, \ldots \, \, ,\, \, N-1,\, \, a_{N} =\varepsilon _{N-1}$. \\

The standard linearization scheme for \eqref{3} looks as following.
\begin{equation} \label{4} 
\, \begin{array}{c} {x_{n+1} =a_{1} f\left(x_{n} \right)+a_{2} f\left(x_{n-2} \right)+\, ...\, \, +a_{N} f\left(x_{n-2(N-1)} \right)} \\ {x_{n+2} =a_{1} f\left(x_{n+1} \right)+a_{2} f\left(x_{n-1} \right)+\, ...\, \, +a_{N} f\left(x_{n-2(N-1)+1} \right)} \end{array} 
\end{equation} 
The solution to \eqref{4} can be written as 
\begin{equation} \label{5} 
\, \begin{array}{c} {x_{2n} =\eta _{1} +u_{n} } \\ {x_{2n+1} =\eta _{2} +v_{n} } \end{array}.                                                                          
\end{equation} 
So, if we introduce a vector 
$\, Y_{n} =\left(\begin{array}{c} {x_{2n} } \\ {x_{2n+1} } \end{array}\right)$, the cycle corresponds to the equilibrium 
$\, \left(\begin{array}{c} {\eta _{1} } \\ {\eta _{2} } \end{array}\right)$ w.r.t. this vector.\\

Let put a solution \eqref{5} to the system \eqref{4} and assume that the  increments $u_{n} , v_{n}$ are small.\\

And let $n=2m,$ then $x_{n+1} =x_{2m+1} =\eta _{2} +\, v_{m} $, $x_{n+2} =x_{2m+2} =\eta _{1} +\, u_{m+1} $ and
\[\, \begin{array}{c} {\eta _{2} +v_{m} =a_{1} (f\left(\eta _{1} \right)+f'\left(\eta _{1} \right)u_{m} )+\, ...\, \, +a_{N} (f\left(\eta _{1} \right)+f'\left(\eta _{1} \right)u_{m-N+1} )} 
\\ 
{\eta _{1} +u_{m+1} =a_{1} (f\left(\eta _{2} \right)+f'\left(\eta _{2} \right)v_{m} )+\, ...\, \, +a_{N} (f\left(\eta _{2} \right)+f'\left(\eta _{2} \right)\delta _{m-N+1} )} \end{array}.\] 
Since $\, \, \eta _{1} =f\left(\eta _{2} \right), \eta _{2} =f\left(\eta _{1} \right)$ then
\begin{equation} \label{GrindEQ__9_} 
\, \begin{array}{c} {v_{m} =f'\left(\eta _{1} \right)(a_{1} u_{m} +a_{2} u_{m-1} +\, ...\, \, +a_{N} u_{m-N+1} )} \\ {u_{m+1} =f'\left(\eta _{2} \right)(a_{1} v_{m} +a_{2} v_{m-1} +\, ...\, \, +a_{N} v_{m-N+1} )} \end{array}.                                          
\end{equation} 
The system \eqref{GrindEQ__9_} is linear therefore it solutions can be written as
$$
\begin{array}{c} {v_{m} =c_{1} \lambda ^{m} } \\ {u_{m} =c_{2} \lambda ^{m} } \end{array},
$$ 
which after substitution to \eqref{GrindEQ__9_} gives
\[\, \begin{array}{c} {c_{1} =f'\left(\eta _{1} \right)(a_{1} +a_{2} \lambda ^{-1} +\, ...\, \, +a_{N} \lambda ^{-N+1} )\cdot c_{2} } \\ {\lambda \cdot c_{2} =f'\left(\eta _{2} \right)(a_{1} +a_{2} \lambda ^{-1} +\, ...\, \, +a_{N} \lambda ^{-N+1} )\cdot c_{1} } .\end{array}
\]
That might be written in a vector form as $A\cdot C=O$ where $A, C$ and $O$ are the following matrices
\begin{equation} \label{GrindEQ__10_} 
A=\left( \begin{array}{cc} {1} & {-f'\left(\eta _{1} \right)p(\lambda)} \\ {-f'\left(\eta _{2} \right)p(\lambda)} & {\lambda } \end{array}\right),\quad
C=\left(\begin{array}{c} {c_{1} } \\ {c_{2} } \end{array}\right),\quad
O=\left(\begin{array}{c} {0} \\ {0} \end{array}\right)          
\end{equation} 
and $p(\lambda)=a_{1} +a_{2} \lambda ^{-1} +\, ...\, \, +a_{N} \lambda ^{-N+1}.$

Since $c_{1}, c_{2} $ cannot be both equal to zero then the determinant of the matrix $A$ should be zero, i.e.
\begin{equation} \label{8} 
\lambda -f'\left(\eta _{1} \right)f'\left(\eta _{2} \right)(a_{1} +a_{2} \lambda ^{-1} +\, ...\, \, +a_{N} \lambda ^{-N+1} )^{2} =0.                                              
\end{equation} 
The equation \eqref{8} is a characteristic equation for the system \eqref{3}. If all the roots are in the unit disc then the cycle 
$\left(\eta _{1} ,\, \, \eta _{2} \right)$ is asymptotically stable.\\

The characteristic polynomial can be reduced to the standard form by the substitution $\nu =\lambda ^{1/2}$:
\begin{equation} \label{GrindEQ__12_} 
f(\nu )=\nu ^{2N-1} +k\, (a_{1} \nu ^{2N-2} +a_{2} \nu ^{2N-4} ...+a_{N} ),                                           
\end{equation} 
where $k=\pm i\sqrt{\left|\mu \right|} $ (here we took into account that the cycle multiplier $\mu =f'\left(\eta _{1} \right)f'\left(\eta _{2} \right)$ is negative). \\

For a sufficiently small $\left| \mu  \right|$ all the roots of the characteristic equation (5) are withing the standard unit circle in the complex plane. With $\left| \mu  \right|$ is increasing   certain roots could appear on a circle, then proceed to the exterior of a unit disc, and/or come back inside the circle etc. For the problem of chaos stabilizing the first value of $\left| \mu  \right|$ that allows roots to appear on a circle and to  get of to the exterior plays a crucial role. On a circle the following equation is valid
\begin{equation} \label{GrindEQ__13_} 
\frac{1}{\mu } = {e^{ - i\omega }}\,{\left( {p\left( {{e^{i\omega }}} \right)} \right)^2}\quad
\end{equation}                                               
or                                                           
$$
\frac{1}{k} = ({a_1}{e^{ - i\omega }} + {a_2}{e^{ - i3\omega}} + \,\,\,...\,\,\, + 
{a_N}{e^{ - i(2N -1)\omega}}).
$$                                  
For positive $\mu$ because 
$$
\sum\limits_{j = 1}^N {{a_j} = 1}
$$
if the system is stable then $\mu$ does not exit one.\\

However, if $\mu$ is negative then choosing ${a_1},\,\, \ldots \,\,,\,{a_N}$ appropriately the limited value of $\left| \mu  \right|$ can be more then 1. Denote that value by $\mu ^ *.$ In this notations if $\mu  \in \left( { - {\mu ^ *},\,\,0} \right)$ then all equations of the family (5) have the roots inside the standard unit disc, and for any positive $\varepsilon$ if $\mu  \in \left( { - {\mu ^ * } - \varepsilon ,\,\,0}\right)$ then there are equations in (5) whose roots are outside the standard unit disc.\\

{\bf Problem 1.}
Find 
$$
\mathop {\sup }\limits_{{a_1},\, \ldots
\,,\,{a_N}} \{ {\mu ^ * }\left( {{a_1},\,\, \ldots \,\,,\,{a_N}}
\right)\}  = \bar \mu. 
$$

{\bf Problem 2.} Find the existence of the control of the form (2) such that each equation in the system (3) has stable cycles, find
$$
\mathop {\min }\limits_{{a_1},\, \ldots \,,\,{a_N}} \{ {N^
* }\left( {{a_1},\,\, \ldots \,\,,\,{a_N}} \right)\}  = 
\bar N$$ 
and find the optimal coefficients.\\

The problems 1 and 2  are dual to each other.\\

{\bf Preliminary results.}
To construct an output function change \eqref{GrindEQ__13_}  to the form
$$\frac{1}{\mu } = {\left( {{a_1}{e^{ - it}} + {a_2}{e^{ - i3t}} + ... + {a_N}{e^{ - i(2N + 1)t}}} \right)^2},$$
where $t = \frac{\omega }2.$

Denote by $\Omega(t)$ the set
$$
\Omega(t) :=
\left\{ {\rm Re} 
\left( 
\sum\limits_{j = 1}^N a_j e^{-i(2j+1)t}
\right)^2:
{\rm Im}
\left( 
\sum\limits_{j = 1}^N a_j e^{-i(2j+1)t}  
\right)^2 = 0 
\right\}
$$
and let
\begin{equation} \label{11}
I\left(N \right) = 
\mathop {\sup }\limits_{\sum\limits_{j = 1}^N a_j =1} 
\left[ 
\min_{t \in \left[0,2\,\pi\right]} \Omega(t)
\right]
\end{equation}
Because $\cos $ is an even and $\sin$ is an odd function the minimum in \eqref{11} 
could be taken along $\left[ {0,\,\pi }\right].$\\

{\bf Lemma 1}. {\sl The quantity 
$$
\min_{t \in \left[0,\pi\right]} \Omega(t)
$$ is negative.}\\

{\bf Proof}. For  the function 
$$
F(z) = \sum\limits_{j = 1}^N {a_j}{z^{2j+1}}
$$
of a complex variable $z$ the point ${z_0} = 0$ is a root. Because of continuity of a polynimial roots
on the coefficients the function 
$$
{F_\varepsilon }(z) = \varepsilon  + \sum\limits_{j =
1}^N {a_j}{z^{2j+1}}
$$ 
does have zero point ${z_\varepsilon }$ whose absolute value for sufficiently small $\varepsilon$ is less then one.

Let 
$$
\min_{t \in \left[0,\pi\right]} \Omega(t) \ge 0.
$$
Then for $\varepsilon>0$ the function graph 
$$
x + iy = {F_\varepsilon}({e^{it}})
$$
does not intersect the real negative half axis on ОХУ plane and 
doesn't pass through the origin, i.e. doesn't surround any zero 
for $t \in [\,0,\,2\pi \,]$. Therefore, by the argument principle
the function ${F_\varepsilon }(z)$ doesn't have zeros inside the unit disc which is wrong.
Lemma 1 is proved.\\

{\bf Lemma 2} 
$$
\bar \mu  =  - \frac{1}{{I\left(N\right)}}.
$$

The proof follows from the definition of $I\left(N\right).$ \\

Let 
$$
J\left(N \right) = 
\mathop {\inf} \limits_{\sum\limits_{j = 1}^N {{a_j}}  = 1} 
\left[ {\mathop{\max }\limits_{t \in \left[ {0,\,\frac{\pi }{2}} \right]} 
\left\{
{\left| {\sum\limits_{j = 1}^N a_j e^{-i(2j - 1)t}}
\right|:\,\arg \left( {\sum\limits_{j = 1}^N a_j e^{-i(2j - 1)t} } \right) 
= \frac{\pi }{2}\,} \right\}} 
\right].
$$
It is easy to show that the set
$$
\left\{
{\left| {\sum\limits_{j = 1}^N a_j e^{-i(2j - 1)t}}
\right|:\,\arg \left( {\sum\limits_{j = 1}^N a_j e^{-i(2j - 1)t} } \right) 
= \frac{\pi }{2}\,}
\right\}
$$ 
is not empty, and moreover   
$$
I\left(N\right) =  - {\left( {J\left(N\right)} \right)^2}.
$$

 To calculate the values ​​of  
$\bar \mu = {\mu ^ * }\left( {a_1^0,\, \ldots \,,\,a_N^0} \right),$
$a_1^0,\,\ldots \,,\,a_N^0$, ${N^ * }$ and build the optimal control of the form (2), which will stabilize a cycle, we will use the procedure developed in [2]. In that work, problems 1 and 2 are completely solved for an equlibrium point.\\

One can write 
$$
J\left(N\right) = \mathop {\inf }\limits_{\sum\limits_{j = 1}^N {{a_j}}  = 1} 
\left[ {\mathop {\max }\limits_{t \in \left[ {0,\frac{\pi }{2}} \right]} \left\{ {|S(t)|:C(t) = 0} 
\right\}} \right]
$$
where
$$
C(t)=\sum_{j=1}^Na_j\cos(2j-1)t,\quad S(t)=\sum_{j=1}^Na_j\sin(2j-1)t.
$$
For $t \in \left( {0,\,\frac{\pi }{2}} \right)$ we have
$$
C(t)=\frac1{2\sin t}\sum_{j=1}^N \alpha_i\sin2jt
$$
where  ${\alpha _j} = {a_j} - {a_{j + 1}},\,\,j = 1,\, \ldots
\,,\,N$ (and we assume that ${a_{N + 1}} = 0.$)\\

In [2]  Lemmas 2-5 which containe some properties of the conjugate trigonometric polynomials
were were formulated and proved. Since the trigonometric polynomials  $2\sin t C(t)$ and $\sum\limits_{j = 1}^N {{\alpha _j} - 2\sin t\,S(t)} $ are conjugate those lemmas implies direct analogous properties of trigonometric polynomials $C(t)$ and $S(t).$ Let formulate them.
\\

{\bf Lemma 3.} {\sl Let $C({t_1}) = 0,\,\,{t_1} \in \left( {0,\,\frac{\pi}{2}} \right).$ Then trigonometric polynomials $C(t)$ and $S(t)$ are represented uniquely in the form
$$
C(t) = \frac{1}{{2\sin t}}(\cos 2t - \cos 2{t_1}) 
\cdot \sum\limits_{j = 1}^{N - 1} {\alpha _j^{(1)}\sin 2jt},
$$ 
$$
S(t) = \frac{1}{{2\sin t}}\left( {\sum\limits_{j = 1}^N {{\alpha _j}}  + \frac{{\alpha _1^{(1)}}}{2} - (\cos 2t - \cos 2{t_1}) \cdot 
\sum\limits_{j = 1}^{N-1} {\alpha _j^{(1)}\cos 2jt\,\,} } \right).
$$
}
\\

{\bf Lemma 4.} {\sl 
If $C({t_1}) = 0,\,\,{t_1} \in \left( {0,\,\frac{\pi
}{2}} \right)$ then 
$$
S({t_1}) = \frac{1}{{2\sin {t_1}}}\left(
{\sum\limits_{j = 1}^N {{\alpha _j}}  + \frac{{\alpha
_1^{(1)}}}{2}} \right),
$$ 
where the values ${\alpha_1},\,...\,,\,{\alpha _N},\;\alpha _1^{(1)}$ are uniquely determined by the coefficients ${a_1},\,...\,,\,{a_N}.$
}
\\

{\bf Lemma 5.} {\sl If
$$
C({t_1}) = 0,\,\,{t_1} \in \left( {0,\,\frac{\pi }{2}} \right),
$$
then
$$
S(\frac{\pi }{2}) = \frac{1}{2}\left( {\sum\limits_{j = 1}^N {{\alpha _j}}  + \frac{{\alpha _1^{(1)}}}{2} + (1 + \cos 2{t_1}) \cdot \sum\limits_{j = 1}^{N - 1} {{{( - 1)}^j} \cdot \alpha _j^{(1)}\,\,} } \right).
$$
}
\\

{\bf Lemma 6.} {\sl Let 
$$
C({t_1}) = \,\, \ldots \,\, = C({t_m}) =
0,\;{t_j} \in \left( {0,\,\frac{\pi }{2}} \right),\,\,j = 1,\ldots m\;(m < N).
$$ 
Then $C(t)$ can be uniquely represented as
$$
C(t) = \frac{1}{{2\sin t}}\prod\limits_{j = 1}^m {(\cos 2t - \cos 2{t_j})}  
\cdot \sum\limits_{j = 1}^{N - m} {\alpha _j^{(m)}} \sin 2jt,
$$
where the values  $\alpha_j^{(m)},\,j = 1,\ldots ,N - m$ are uniquely determined by the coefficients ${a_1},...,{a_N}.$
} \\

{\bf Lemma 7.} {\sl Let $C(t_1)=\dots=C(t_m)=0, S(t_1)=\dots=S(t_m)$ , where ${t_1},\,\, \ldots
\,,{t_m}$- are distinct and belong to the interval 
$\left({0,\,\frac{\pi }{2}} \right).$ Then $S(t)$ is uniquely represented as
$$
S(t) = \frac{1}{{2\sin t}}\left( {\sum\limits_{j = 1}^N {{\alpha _j}}  + \frac{{\alpha _m^{(m)}}}{{{2^m}}} - \prod\limits_{j = 1}^m {(\cos 2t - \cos 2{t_j})}  \cdot \sum\limits_{j = m}^{N - m} {\alpha _j^{(m)}} \cos 2jt} \right),
$$
where the values $\alpha _j^{(m)},\,j = m,\, \ldots \,,\,N - m$
are uniquely determined by the coefficients ${a_1},...,{a_N}.$
}
\\

{\bf Lemma 8.} {\sl If $C(t_1)=\dots=C(t_m)=0,\; S({t_1}) = \,\, \ldots \, = S({t_m})\,$, 
$m < N,$ then 
$$
S({t_1}) = \,\, \ldots \, = S({t_m}) = \frac{1}{{2\sin {t_j}}}\left( {\sum\limits_{j = 1}^N {{\alpha _j}}  + \frac{{\alpha _m^{(m)}}}{{{2^m}}}} \right)\,
$$ where the value  $\alpha _m^{(m)}$ are uniquely determined by the coefficients 
${a_1},...,{a_N}.$
}
\\

{\bf Lemma 9.} {\sl If
$S({t_1}) = \,\, \ldots \, = S({t_m}),$ $m < N$, then 
$$
S(\frac{\pi }{2}) = \frac{1}{2}\left( {\sum\limits_{j = 1}^N {{\alpha _j}}  + \frac{{\alpha _m^{(m)}}}{{{2^m}}} + \prod\limits_{j = 1}^m {(1 + \cos 2{t_j})}  \cdot \sum\limits_{j = m}^{N - m} {{{( - 1)}^j}\alpha _j^{(m)}} } \right).
$$
}
\\

{\bf 5. Main results.} 

{\bf Theorem.} {\sl The following extremal values are found  
$$  
- I\left( N \right) = {\left( {J\left(N \right)} \right)^2} = \frac{1}{{{N^2}}}.
$$
}

{\bf Proof}. Let $C(t),\,\,S(t)$ be a pair of conjugate trigonometric polynomials
$$
C(t) = \sum\limits_{j = 1}^N {{a_j}\cos (2j - 1)t,\quad } S(t) = \sum\limits_{j = 1}^N {{a_j}\sin (2j - 1)t,} 
$$
normalized by the condition 
$$\sum\limits_{j = 1}^N {{a_j} = 1}.$$
We consider the extremal problem
$$
\rho  = \mathop {\inf }\limits_{\sum\limits_{j = 1}^N {{a_j}}  = 1} \left[ {\mathop {\max }\limits_{t \in \left[ {0,\frac{\pi }{2}} \right]} \left\{ {|S(t)|:C(t) = 0\,} \right\}} \right].
$$
Let  $T$ denote denote the set of sign changes of the polynomial $C(t)$ on the interval 
$(\,0,\,\frac{\pi }{2}\,)$, and let 
$$
{\rho _1} = \mathop {\inf }\limits_{\sum\limits_{j = 1}^N {{a_j}}  = 1} \left[ {\mathop {\max }\limits_{} \left\{ {|S(t)|:t \in T \cup \left\{ {\frac{\pi }{2}} \right\}\,} \right\}} \right].
$$
The same way as it is done in the proof of Theorem 1 [2], it can be shown, that infimum for ${\rho _1}$ can be reached, therefore $\rho  \ge {\rho _1}$ and there exists a pair of optimal polynomials. So let 
$\left\{ {{C^0}(t),{S^0}(t)}\right\}$  be an optimal pair and  ${T^0}$ be the set of sign changes of the polynomial $C(t)$ on the interval 
$(\,0,\,\frac{\pi }{2}\,).$ Then
$$
{\rho _1} = \mathop {\max }\limits_{} \left\{ {\,\,{S^0}(t):t \in
{T^0} \cup \left\{ {\frac{\pi }{2}} \right\}\,} \right\}.
$$
Let show, that ${T^0} = \emptyset $ thus 
$$
{\rho _1}
= \left| {{S^0}(\frac{\pi }{2})} \right|,
$$
and ${C^0}(t) \ge 0$ on $(\,0,\,\frac{\pi }{2}).$\\

Suppose that for the optimal polynomial ${C^0}(t)$ the set $T=\{t_1,\dots,t_q\},$ where
$0 \le q \le n - 1$ , is not empty. And let $\max \left\{ {\,\,{S^0}({t_1}),...,{S^0}({t_q})\,}
\right\} = {S^0}({t_1})$, and 
$$
{S^0}({t_1}) = {S^0}({t_j}),j = 1,...,m\;(1 \le m \le q),\,{S^0}({t_1}) > {S^0}({t_j}),\,\,j = m + 1,\,\, \ldots \,,\,q.
$$
Two cases are possible: 
$$
{S^0}({t_1}) > {S^0}(\frac{\pi }{2}),\,{S^0}({t_1})\le {S^0}(\frac{\pi }{2}).
$$

{\bf Case 1.} Accordingly to lemmas 6,7 the trigonometric polynomials ${S^0}(t),\,\,{C^0}(t)$ have the form
$$
{C^0}(t) = \frac{1}{{2\sin t}}\prod\limits_{j = 1}^m {(\cos 2t - \cos 2{t_j})}  \cdot \sum\limits_{j = m}^{N - m} {\alpha _j^{(m)}} \sin 2jt,
$$
$$
{S^0}(t) = \frac{1}{{2\sin t}}\left( {\sum\limits_{j = 1}^N {{\alpha _j}}  + \frac{{\alpha _m^{(m)}}}{{{2^m}}} - \prod\limits_{j = 1}^m {(\cos 2t - \cos 2{t_j})}  \cdot \sum\limits_{j = m}^{N - m} 
{\alpha _j^{(m)}} \cos 2jt} \right).
$$
Since 
$$
{S^0}({t_1}) = \frac{1}{{2\sin {t_1}}}\left(
{\sum\limits_{j = 1}^N {{\alpha _j}}  + \frac{{\alpha_m^{(m)}}}{{{2^m}}}} \right),
$$
then
$$
\sum\limits_{j = 1}^N{{\alpha _j}}  + \frac{{\alpha _m^{(m)}}}{{{2^m}}} > 0.
$$
Since ${C^0}(0) = 1$ then
$$
\prod\limits_{j = 1}^m{(1 - \cos 2{t_j})}  \cdot \sum\limits_{j = m}^{n - m} 
{j\alpha_j^{(m)}}  = 1,
$$   
$$
\sum\limits_{j = m}^{n - m} {j\alpha_j^{(m)}}  > 0.
$$

 Let’s construct auxiliary polynomials
$$
C({\theta _1},\, \ldots \,,\,{\theta _m};\,t) = N({\theta _1},\, \ldots \,,\,{\theta _m}) \cdot \frac{1}{{2\sin t}}\prod\limits_{j = 1}^m {(\cos 2t - \cos 2{\theta _j})} \sum\limits_{j = m}^{n - m} {\alpha _j^{(m)}\sin 2jt}, 
$$
$$
S({\theta _1},\, \ldots \,,\,{\theta _m};\,t) = N({\theta _1},\, \ldots \,,\,{\theta _m}) \cdot \frac{1}{{2\sin t}}\cdot
$$
$$
\left( {\sum\limits_{j = 1}^N {{\alpha _j}}  + \frac{{\alpha _m^{(m)}}}{{{2^m}}} - \prod\limits_{j = 1}^m {(\cos 2t - \cos 2{\theta _j})}  \cdot \sum\limits_{j = m}^{N - m} {\alpha _j^{(m)}} \cos 2jt} \right),
$$
where the normalizing factor $N({\theta _1},\, \ldots \,,\,{\theta_m})$ makes sum of polynomial coefficients $S({\theta _1},\, \ldots \,,\,{\theta _m};\,t),$
$C({\theta _1},\, \ldots \,,\,{\theta _m};\,t)$ to be equal one. For the polynomial
$C({\theta _1},\, \ldots \,,\,{\theta _m};\,t)$ the set of sign changes 
${T_\theta } = \left\{ {{\theta _1},\, \ldots
\,,\,{\theta _m},\,{t_{m + 1}},\, \ldots \,,\,{t_q}} \right\}.$\\

It is clear, that $S({t_1},\, \ldots \,,\,{t_m};\,t) \equiv
{S^0}(t),\,\,C({t_1},\, \ldots \,,\,{t_m};\,t) \equiv {C^0}(t).$
The factor $N({\theta _1},\, \ldots \,,\,{\theta _m})$  is determined by the condition 
$C({\theta _1},\, \ldots \,,\,{\theta_m};\,0) = 1$, i.e. 
$$
N({\theta _1},\, \ldots \,,\,{\theta _m}) = \frac{1}{{\prod\limits_{j = 1}^m {(1 - \cos 2{\theta _j})} \sum\limits_{j = m}^{n - m} {j\alpha _j^{(m)}} }}.
$$

Polynomials $S({\theta _1},\, \ldots \,,\,{\theta _m};\,t)$  finally could be defined as
$$
C({\theta _1},\, \ldots \,,\,{\theta _m};\,t) = \frac{1}{{\prod\limits_{j = 1}^m {(1 - \cos 2{\theta _j})} \sum\limits_{j = m}^{n - m} {j\alpha _j^{(m)}} }} \cdot
$$
$$ 
\frac{1}{{2\sin t}}\prod\limits_{j = 1}^m {(\cos 2t - \cos 2{\theta _j})} \sum\limits_{j = m}^{n - m} {\alpha _j^{(m)}\sin 2jt}. 
$$
$$
S({\theta _1},\, \ldots \,,\,{\theta _m};\,t)
= \frac{1}{{\prod\limits_{j = 1}^m {(1 - \cos 2{\theta _j})} \sum\limits_{j = m}^{n - m} {j\alpha _j^{(m)}} }} \cdot
$$ 
$$
\frac{1}{{2\sin t}}\left( {\sum\limits_{j = 1}^N {{\alpha _j}} +
\frac{{\alpha _m^{(m)}}}{{{2^m}}} - \prod\limits_{j = 1}^m {(\cos 2t - \cos 2{\theta _j})}  \cdot \sum\limits_{j = m}^{N - m} {\alpha _j^{(m)}} \cos 2jt} \right).
$$
Let show that the value of ${\rho _1}$ for the pair
$$
\left\{ {S({\theta
_1},\, \ldots \,,\,{\theta _m};\,t),\,C({\theta _1},\, \ldots
\,,\,{\theta _m};\,t)} \right\}
$$
is less than for the pair $\left\{{{S^0}(\,t),\,{C^0}(t)} \right\},$ i.e. the pair
$\left\{{{S^0}(\,t),\,{C^0}(t)} \right\}$ cannot be optimal. 
From Lemma 8, we obtain
$$
S({\theta _1},\, \ldots \,,\,{\theta _m};\,{\theta _j}) = \,\,\frac{{\left( {\sum\limits_{j = 1}^N {{\alpha _j}}  + \frac{{\alpha _m^{(m)}}}{{{2^m}}}} \right)}}{{\prod\limits_{j = 1}^m {(1 - \cos 2{\theta _j})} \sum\limits_{j = m}^{n - m} {j\alpha _j^{(m)}} }} \cdot \frac{1}{{2\sin {\theta _j}}}\,,\,j = 1,\, \ldots \,,\,m.
$$	
Since
$$
\left( {\sum\limits_{j = 1}^N {{\alpha _j}}  +
\frac{{\alpha _m^{(m)}}}{{{2^m}}}} \right) > 0
$$ 
and
$$
\sum\limits_{j = m}^{n - m} {j\alpha _j^{(m)}}  > 0,
$$
so all the values 
$$
S({\theta _1},\, \ldots \,,\,{\theta _m};\,{\theta
_j}),\; j = 1,\, \ldots \,,\,m,
$$ 
decrease by the parameters ${\theta _1},\, \ldots \,,\,{\theta _m}.$ 
Let 
$$
0 < {\theta_j} - {t_j} < \varepsilon ,\,j = 1,\,\, \ldots \,,\,m.
$$
Then, by the continuity of trigonometric polynomials  on $t$ and on all coefficients we get the inequalities  
$$
S({\theta _1},\, \ldots
\,,\,{\theta _m};\,{\theta _j}) < {S^0}({t_j}),\,j = 1,\, \ldots
\,,\,m\,
$$ 
$$
\left| {S({\theta _1},\, \ldots \,,\,{\theta
_m};\,{t_j}) - {S^0}({t_j})} \right| < \delta ,\,j = m + 1,\,\,
\ldots \,\,,\,q,
$$
$$
\left| {S({\theta _1},\, \ldots \,,\,{\theta _m};\,\frac{\pi }{2}) - {S^0}(\frac{\pi }{2})} \right| < \delta 
$$ 
for an arbitrarily small $\delta$ with an appropriate choice $\varepsilon.$ 

These inequalities mean that the value
$$
\max \left\{S({\theta _1},\, \ldots \,,\,{\theta _m};\,{\theta _1}),\, \ldots \,,\,S({\theta _1},\, \ldots \,,\,{\theta _m};\,{\theta _m}),\right.
$$
$$
\left.
S({\theta _1},\, \ldots \,,\,{\theta _m};\,{t_{m + 1}}),\,\, \ldots \,,\,S({\theta _1},\, \ldots \,,\,{\theta _m};\,{t_q})\,,\, \ldots \,,\,S({\theta _1},\, \ldots \,,\,{\theta _m};\,{\pi }/{2}) \right\}
$$
is less then 
$$
\max \left\{
{\,\,{S^0}({t_1}),...,{S^0}({t_q}),{S^0}(\frac{\pi }{2})\,}
\right\}
$$ at least for sufficiently small positive 
$${\theta _j} - {t_j},\,j = 1,\,\, \ldots \,,\,m,$$ i.e. the pair $\left\{ {{S^0}(\,t),\,{C^0}(t)} \right\}$ is not optimal.\\

{\bf Case 2.} From Lemma 9, we obtain
$$
{S^0}(\frac{\pi }{2}) = \frac{1}{2}\left( {\sum\limits_{j = 1}^N {{\alpha _j}}  + \frac{{\alpha _m^{(m)}}}{2} + \prod\limits_{j = 1}^m {(1 + \cos 2{t_j})}  \cdot \sum\limits_{j = m}^{N - 1} {{{( - 1)}^j} \cdot \alpha _j^{(1)}\,\,} } \right),
$$
$$
S({\theta _1},\, \ldots \,,\,{\theta _m};\,\frac{\pi }{2}) = 
\frac{1}{ {2\prod\limits_{j = 1}^m {(1 - \cos 2{\theta _j})} 	
\sum\limits_{j = m}^{n - m} {j\alpha _j^{(m)}} }}\cdot
$$
$$
\left( {\sum\limits_{j = 1}^N {{\alpha _j}}  + \frac{{\alpha _m^{(m)}}}{2} + \prod\limits_{j = 1}^m {(1 + \cos 2{\theta _j})}  \cdot \sum\limits_{j = m}^{N - 1} {{{( - 1)}^j} \cdot \alpha _j^{(1)}\,\,} } \right),
$$
and 
$$
S({t_1},\, \ldots \,,\,{t_m};\,\frac{\pi }{2})={S^0}(\frac{\pi }{2}).
$$
Since by assumption 
$$
{S^0}(\frac{\pi }{2}) \ge {S^0}({t_j}) = \frac{1}{{2\sin
{t_j}}}\left( {\sum\limits_{j = 1}^N {{\alpha _j}}  +
\frac{{\alpha _m^{(m)}}}{{{2^m}}}} \right),
$$
then
$$
\prod\limits_{j = 1}^m {(1 + \cos 2{t_j})}  \cdot \sum\limits_{j
= m}^{N - 1} {{{( - 1)}^j} \cdot \alpha _j^{(1)}\,\,}  \ge \left(
{\frac{1}{{\sin {t_j}}} - 1} \right)\left( {\sum\limits_{j = 1}^N
{{\alpha _j}}  + \frac{{\alpha _m^{(m)}}}{{{2^m}}}} \right),
$$
whence
$$
\prod\limits_{j = 1}^m {(1 + \cos 2{t_j})}  \cdot \sum\limits_{j = m}^{N - 1} {{{( - 1)}^j} \cdot \alpha _j^{(1)}\,\,}  \ge 0$$ и  $$\sum\limits_{j = m}^{N - 1} {{{( - 1)}^j} \cdot \alpha _j^{(1)}\,\,}  \ge 0.$$
The values 
$$\sum\limits_{j = 1}^N {{\alpha _j}}  + \frac{{\alpha
_m^{(m)}}}{2} + \prod\limits_{j = 1}^m {(1 + \cos 2{\theta _j})}
\cdot \sum\limits_{j = m}^{N - 1} {{{( - 1)}^j} \cdot \alpha
_j^{(1)}\,\,} 
$$ 
and
$$
\frac{1}{{2\prod\limits_{j = 1}^m {(1 - \cos
2{\theta _j})} \sum\limits_{j = m}^{n - m} {j\alpha _j^{(m)}} }}
$$
decrease for each parameter 
$$
{\theta _1},\, \ldots \,,\,{\theta_m}.
$$
So in this case a pair $\left\{{{S^0}(\,t),\,{C^0}(t)} \right\}$ cannot be optimal. Thus, ${T^0} = \emptyset$.\\

For the polynomial $C(t)$ we have the representation
$$
C(t) = \cos t \cdot ({\gamma _1} + 2{\gamma _2}\cos 2t + .. + 2{\gamma _N}\cos 2(N - 1)t),
$$
where 
$$
{\gamma _s} = \sum\limits_{j = s}^N {{{( - 1)}^{s +
j}}{a_j}},\; s = 1,\,\, \ldots \,,\,\,N.
$$
There is a one-to-one correspondence between ${a_1},\,...\,,\,\,{a_N}$ and ${\gamma _1},\,...\,,\,\,{\gamma_N}$ and
$$
{\gamma _1} + 2\sum\limits_{j = 2}^N {{\gamma _j}}  =
\sum\limits_{j = 1}^n {{a_j} = 1}. 
$$ 
Note that 
$$\left| {S(\frac{\pi}{2})} \right| = \left| {\sum\limits_{j = 1}^N 
{{{( - 1)}^{j +1}}{\gamma _j}} } \right| = \left| {{\gamma _1}} \right|.
$$
Then
$${\rho _1} = \mathop {\min }\limits_{{\gamma _1}, \ldots ,{\gamma _N}} \left\{ {\,\,\left| {{\gamma _1}} \right|:\,\,{\gamma _1} + 2\sum\limits_{j = 2}^N {{\gamma _j}}  = 1,\,\,\frac{{C(t)}}{{\cos t}} \ge 0,\,\,t \in (0,\,\,\frac{\pi }{2})\,} \right\}
$$
and the polynomial  
$$
\frac{{C(\frac{\vartheta }{2})}}{{\cos \frac{\vartheta
}{2}}} = {\gamma _1} + 2{\gamma _2}\cos \vartheta  + \, \ldots \,
+ 2{\gamma _N}\cos (N - 1)\vartheta 
$$ 
is non-negative.\\ 

From the properties of Fej\'er kernel  ([3], 6.7, problem 50) follows that for a non-negative trigonometric polynomial of degree  $N - 1$ the ratio of the value of the polynomial at any point to the mean value does not exceed N, moreover the equality holds at $2k\pi, k\in \mathbb Z$ only.\\

In our case the Fej\'er condition turns to the following one
$$
\frac{C(0)}{\gamma_1}\le N
$$
therefore $\gamma_1^0=1/N$ and ${\rho _1} = {1}/{N}.$\\


Let find $\rho$. To do that we consider a one-parameter family of trigonometric polynomials
$$
{S^\varepsilon }(t) = \sum\limits_{j = 1}^N {a_j^\varepsilon } \sin (2j - 1)t,
\qquad
{C^\varepsilon }(t) = \sum\limits_{j = 1}^N {a_j^\varepsilon } \cos (2j - 1)t,
$$
where 
$$
a_1^\varepsilon  = \frac{{a_1^0 + \varepsilon }}{{1 +
\varepsilon }},\,\,a_j^\varepsilon  = \frac{{a_j^0}}{{1 +
\varepsilon }},\,j = 2,\,\, \ldots \,,\,N.$$
It is clear, that
$$
\sum\limits_{j = 1}^N {a_j^\varepsilon  = } \frac{{a_1^0 +
\varepsilon }}{{1 + \varepsilon }} + \frac{{a_2^0}}{{1 +
\varepsilon }} + ... + \frac{{a_n^0}}{{1 + \varepsilon }} = 1,
$$
and
$$
{S^\varepsilon }(t) = \frac{{{S^0}(t)}}{{1 + \varepsilon }} +
\frac{\varepsilon }{{1 + \varepsilon }}\sin t.
$$
For all $t\in (0,\frac{\pi }{2})$
and $\varepsilon  > 0$  the inequality
${C^\varepsilon }(t) > 0$ is valid. So,
$$
\rho  \le \left| {{S^\varepsilon }(\frac{\pi }{2})} \right| = \frac{{\left| {{S^0}(\frac{\pi }{2})} \right|}}{{1 + \varepsilon }} + \frac{\varepsilon }{{1 + \varepsilon }}.
$$
Let $\varepsilon  \to 0 +$, then in the limit 
$$
\rho  \le \left|{{S^0}(\frac{\pi }{2})} \right| = {\rho _1}.
$$ 
Finally, 
$$
\rho  = {\rho _1} = J\left(N\right) = \frac{1}{N}.
$$ 
The Theorem is proved.
\\

{\bf Corollary 1}. Let pair of conjugate  trigonometric polynomials 
$$
C(t) = \sum\limits_{j = 1}^N {{a_j}\cos (2j - 1)t,\quad } 
S(t) = \sum\limits_{j = 1}^N  {{a_j}\sin (2j - 1)t,}
$$
is normalized by the condition
$$
\sum\limits_{j = 1}^N {{a_j} = 1}.
$$
And let $\tilde J\left(N\right)$ is a solution of the extreme problem
$$
\tilde J\left(N\right)=\mathop {\min }\limits_{\sum\limits_{j =1}^N {{a_j}}  = 1} 
\left[ {\max \left\{ {\,\,S(t):t\in T \cup \left\{
{\frac{\pi }{2}} \right\}\,} \right\}} \right],
$$  
where 
$T$ - is a set of sign changes for  the function $C(t)$ on 
$\left({0,\,\frac{\pi }{2}} \right).$

Then there exists a unique solution 
$$
{C^0}(t) = \sum\limits_{j =  1}^N {a_{_j}^0\cos (2j - 1)t},\qquad
{S^0}(t) = \sum\limits_{j = 1}^N {a_{_j}^0\sin (2j - 1)t,}  
$$ 
where 
$$
a_j^0 = \frac{{2(N - j) + 1}}{{{N^2}}},\;j = 1,\,\,...\,\,,N.
$$  
Moreover, 
$$
\tilde  J\left(N\right) = \frac{1}{N}.
$$
To find the coefficients $a_{_1}^0,\, \ldots \,,\,a_{_N}^0$
we use the representation of the extremal Fej\'er polynomial of double argument and with the average equal to $1/N.$ Its coefficients $\gamma^0_1,\dots,\gamma^0_N$ are define uniquelly and
$$
{C^0}(t) = {\left( {\frac{{\sin Nt}}{{N\sin t}}} \right)^2}\cos t =
$$
$$
\cos t\,\left( {\frac{1}{N} + 2\sum\limits_{j = 2}^N {\frac{{N - j + 1}}{{{N^2}}}\cos 2(j - 1)t} } \right).
$$
Since $a_j^0 = \gamma _j^0 + \gamma _{j + 1}^0,$ $j = 1,...,N$ (we assume that $\gamma _{N + 1}^0=0$), 
than
$$
a_j^0 = \frac{{2(N - j) + 1}}{{{N^2}}},\;j = 1,...,N.
$$

{\bf Corollary 2}.
The optimal coefficients in the control (2) are determined in a unique way:  
$$
{\varepsilon_j} = \sum\limits_{k = j + 1}^N {a_k^0} ,j = 1, \ldots,\,N - 1.
$$

Indeed, the formulas $a_{1} =1-\varepsilon _{1} , a_{j} =\varepsilon _{j-1} -\varepsilon _{j} , j=2,\, \, \ldots \, \, ,\, \, N-1,\, \, a_{N} =\varepsilon _{N-1}$ determine bijection between 
$\varepsilon _{1} ,\dots, \varepsilon _{N-1}$ and $a_1,\dots,a_{N-1}$ and it the formula in Corollary 2 is easy to check.

The graph of the function 
$$
x = {\mathop{\rm Re}\nolimits} \left(
{\sum\limits_{j = 1}^N {a_j^0{e^{ - i(2j - 1)t}}} } \right),\,\,y
= {\mathop{\rm Im}\nolimits} \left( {\sum\limits_{j = 1}^N
{a_j^0{e^{ - i(2j - 1)t}}} } \right)
$$ 
for $N = 4, t\in\left[ {0,\,2\pi } \right]$ is displayed on the Figure 1.\\

The fragment of the above graph for $t \in\left[ {0.7,\,2.4} \right]$ is displayed on the Figure 2.\\

{\bf Corollary 3.} If 
$\left( {{\eta_1},{\eta _2}} \right)$ - is a cycle with a multiplier
${f'_h}\left( {{\eta _1}} \right){f'_h}\left( {{\eta _2}} \right)
\in \left( { - {\mu ^ * },\, - 1} \right),$ then there exists a control of the type (2)
which stabilize the cycle and which is optimal w.r.t. the minimal depth of prehistory in the delayed feedback. For this control ${N^ * } = 2({N_0} - 1),$ where ${N_0}$ is the smallest integer satisfies 
$\sqrt {{\mu ^ * }}  < {N_0}.$\\

{\bf Remark.} If ${\mu ^ * } =\bar \mu,$ then ${N^ * } = \bar N.$\\

{\bf 6. Examples.} For an one-parametric logistic map
$${f_h}\left( x \right) = h \cdot x \cdot (1 - x),\,\,\,0 \le h \le 4,$$                                              
we have
$${f_h}:\,\,\left[ {0,\,\,1} \right] \to \left[ {0,\,\,1} \right].$$
If $h \in \left( {1 + \sqrt 6 ,\,\,4} \right]$  then a cycle
$$
\left( {{\eta _1} = \frac{{1 + h - \sqrt {{h^2} - 2h - 3} }}{{2h}},\,\,{\eta _2} = \frac{{1 + h + \sqrt {{h^2} - 2h - 3} }}{{2h}}} \right)
$$ 
is unstable and a multiplier $\mu  \in \left[ { - 4,\, - 1} \right),$ i.e. ${\mu ^ * } =  - 4$.
Therefore ${N_0} = 3, {N^ * } = 4,$ and the optimal strength coefficients are 
$$
{\varepsilon _1} = a_2^0 + a_3^0 = \frac{4}{9},\quad
{\varepsilon _2} = a_3^0 = \frac{1}{9}.
$$

 Quazistochastic dynamic of the solutions of logistic equation are displayed on the Figure 3, while the solutions closed by the control
$$
u =  - {\varepsilon _1}\left( {{f_h}\left( {{x_n}}
\right) - {f_h}\left( {{x_{n - 2}}} \right)} \right) -
{\varepsilon _2}\left( {{f_h}\left( {{x_n}} \right) - {f_h}\left(
{{x_{n - 4}}} \right)} \right),
$$
which stabilize a cycle are displayed on the Figure 4.\\

{\bf Remark}. If ${N_0} = 2$ then $\left| {{\mu^ * }} \right| = N_0^2.$ 
Hence ${N^ * } = 2,$ 
$$
{\varepsilon_1} = \frac{1}{4},
$$ 
and
$$
u =  - \frac{1}{4}\left( {{f_h}\left( {{x_n}} \right) - {f_h}\left( {{x_{n - 2}}} \right)} \right).
$$                                                (12)
Additional analysis indicates that in this case all solutions of the equation
$$
{x_{n + 1}} = {f_h}\left( {{x_n}} \right) - 
\frac{1}{4}\left( {{f_h}\left( {{x_n}} \right) - {f_h}\left( {{x_{n - 2}}} \right)} \right)
$$
are attracting to some small neighborhood of the periodic orbit evaluated for $h = 4,$ i.e. 
$$
\left( {{\eta _1} = \frac{{5 -\sqrt 5 }}{8},\,\,{\eta _2} = \frac{{5 + \sqrt 5 }}{8}} \right),
$$
however not to the orbit itself, i.e. the control (12) wont be stabilizing to the cycle for $h = 4$ (see Figure 5).\\

Let $h = 3.95$. In this case the control (12) is stabilizing and the Figure 6 displays the dynamic of solution of logistic equations with $h=3.95,$ closed by the control (12).

{\bf 6. Acknowledgement.} The authors would like to thank Alexey Solyanik and Paul Hagelstein for the fruitful discussions, contructive crytics and new ideas suggested.

\newpage
\begin{figure}[h]
{\centering
\includegraphics[width=16cm,height=16cm]{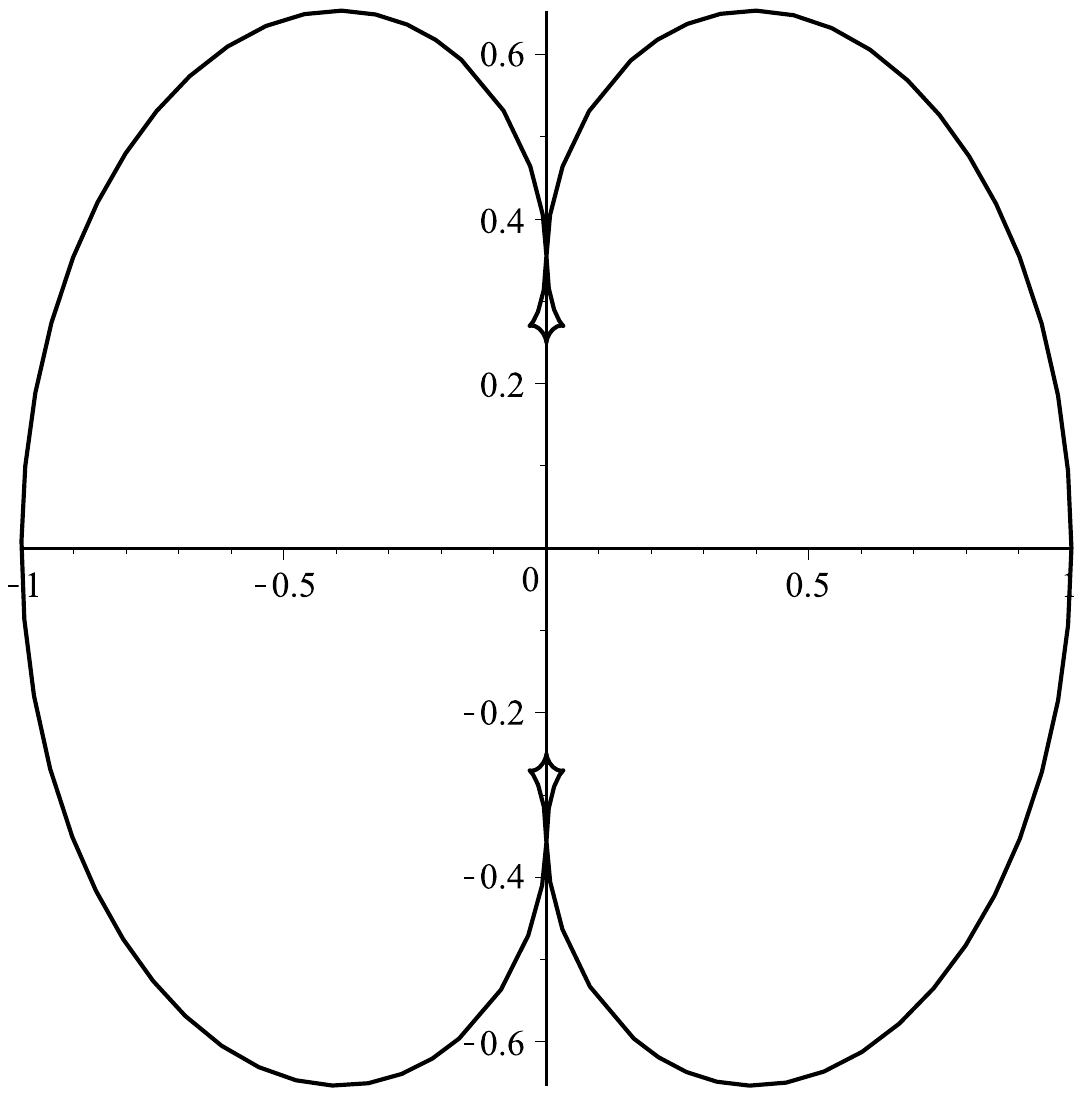}
\caption{}
\label{fig:151106b_measure}
}
\end{figure}

\newpage
\begin{figure}[h]
{\centering
\includegraphics[width=10cm,height=10cm]{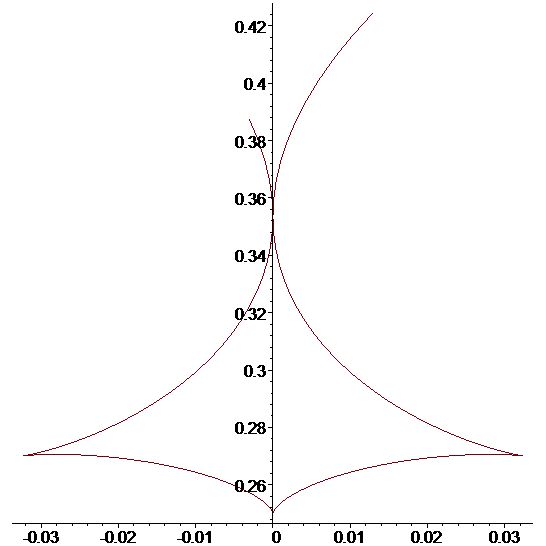}
\caption{}
\label{fig:151106b_measure}
}
\end{figure}

\newpage
\begin{figure}[h]
{\centering
\includegraphics[width=10cm,height=10cm]{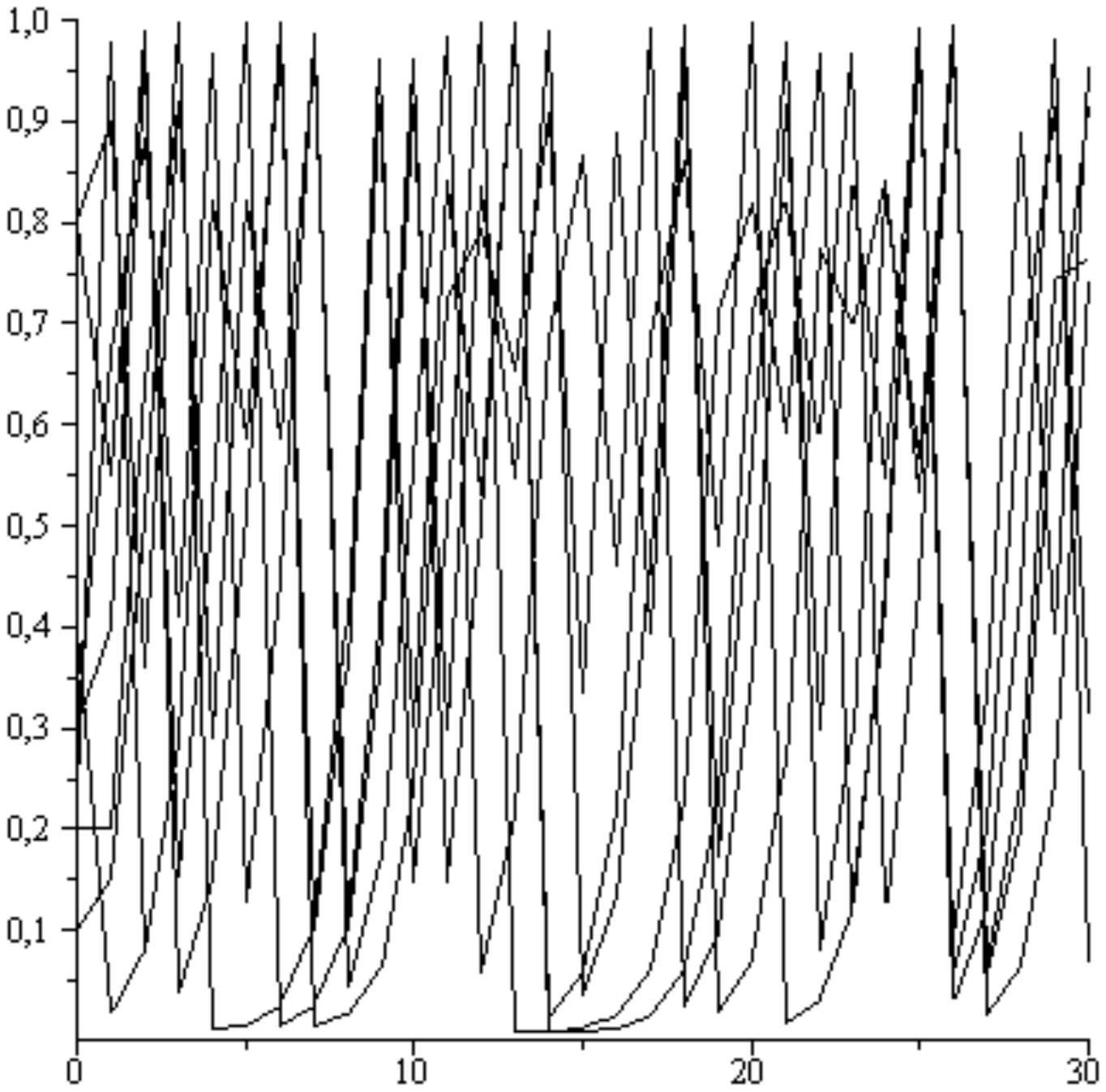}
\caption{}
\label{fig:151106b_measure}
}
\end{figure}

\newpage
\begin{figure}[h]
{\centering
\includegraphics[width=16cm,height=16cm]{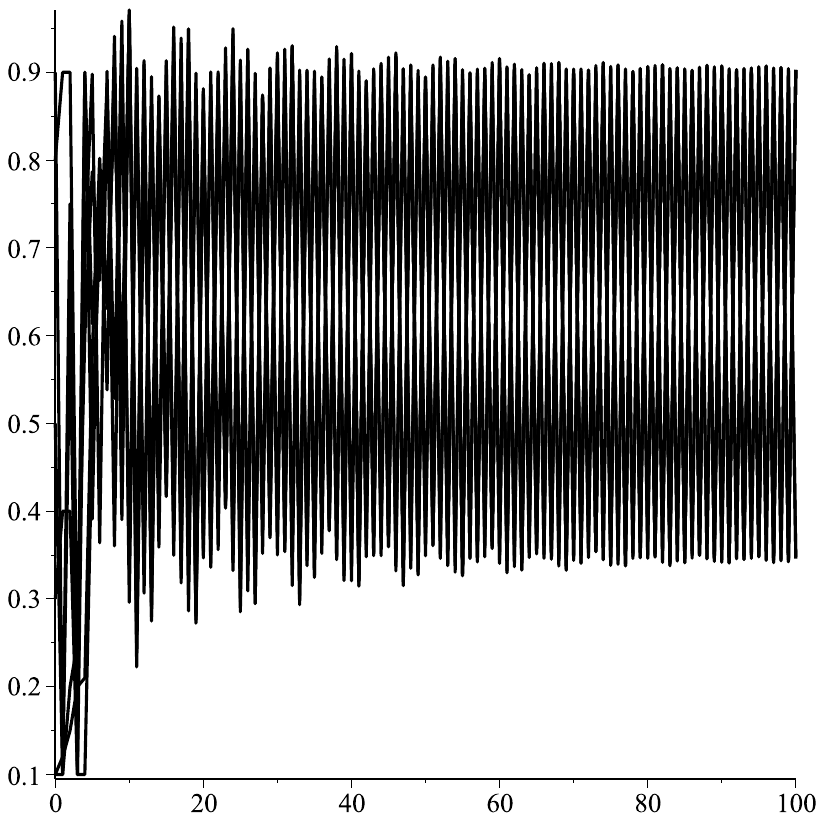}
\caption{}
\label{fig:151106b_measure}
}
\end{figure}

\newpage
\begin{figure}[h]
{\centering
\includegraphics[width=16cm,height=16cm]{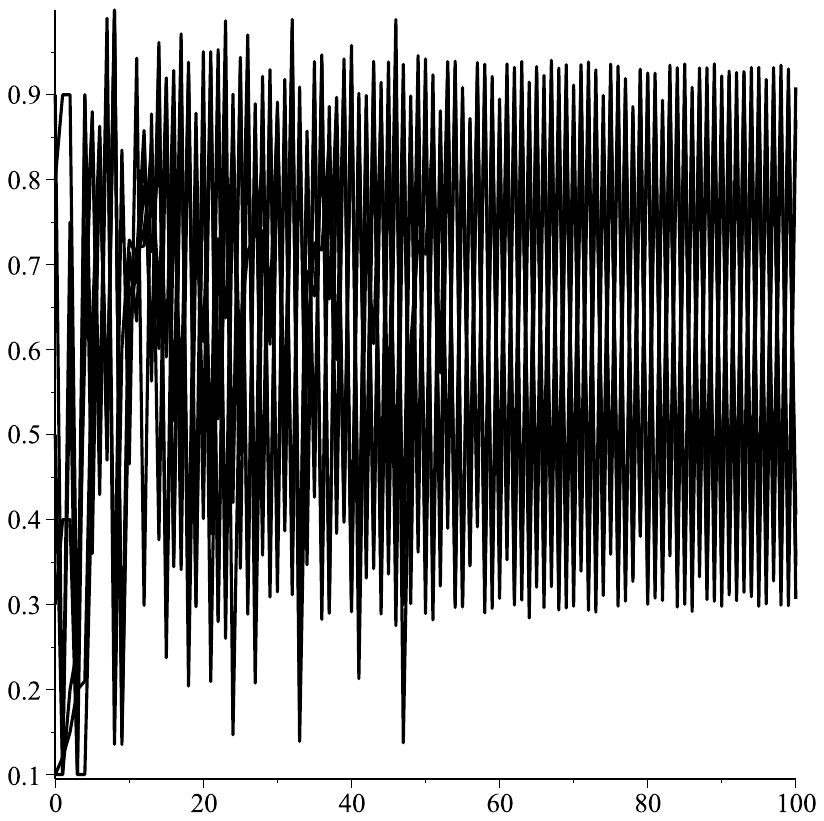}
\caption{}
\label{fig:151106b_measure}
}
\end{figure}

\newpage
\begin{figure}[h]
{\centering
\includegraphics[width=16cm,height=16cm]{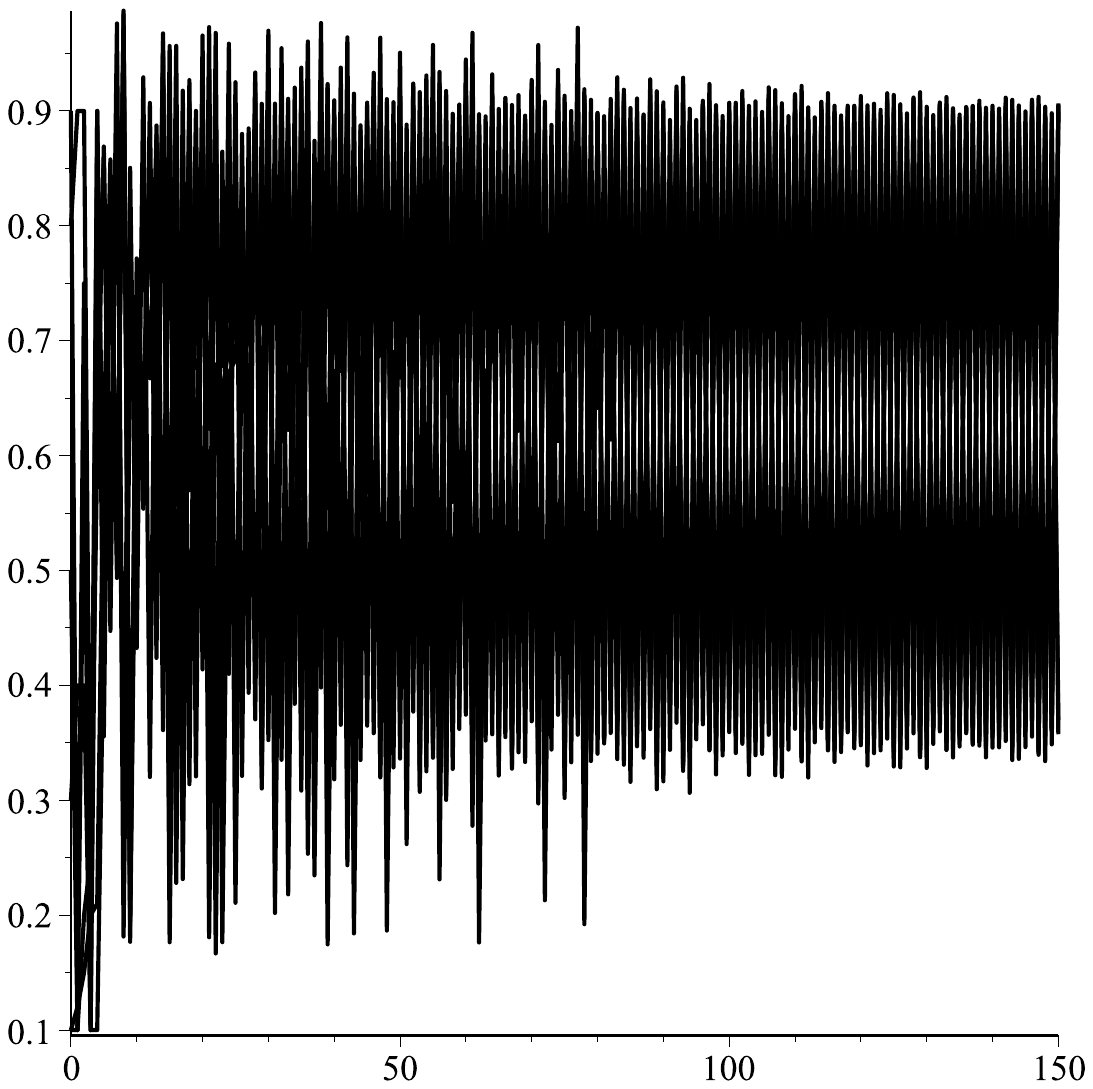}
\caption{}
\label{fig:151106b_measure}
}
\end{figure}

\smallskip

Dmitriy Dmitrishin, Odessa National Polytechnic University, 1 Shevchenko Ave., Odessa 65044, 
Ukraine. E-mail: dmitrishin@opu.ua
      
Anatolii Korenovskyi, Odessa National University, Dvoryanskaya 2, Odessa 65000, Ukraine. E-mail: anakor@paco.net

Alex Stokolos and Anna Khamitova, Georgia Southern University, Statesboro, GA 30458, USA. E-mail: astokolos@georgiasouthern.edu, \\
anna\_khamitova@georgiasouthern.edu

\end{document}